\documentclass{amsart}%%%{article}
\usepackage{amsfonts}
\usepackage{amsmath}
\usepackage{tikz}
\usetikzlibrary{calc}
\usepackage{amssymb}
\usepackage{amstext}
\usepackage{tikz}
\usetikzlibrary{calc}
\usetikzlibrary{decorations.markings}

\newtheorem{theorem}{Theorem}[section]
\newenvironment{PrfFact}{{\bf Proof }}{{\hfill{$\blacksquare$}}}

\theoremstyle{remark}
\newtheorem{remark}{Remark }

\begin{document}

\title{Arc-Disjoint Cycles and Feedback Arc Sets}
\author{Jan~Florek}
\email{jan.florek@ue.wroc.pl}

\address{Institute of Mathematics\\
University of Economics \\
ul. Komandorska 118/120
\\ 53--345 Wroc{\l}aw, Poland
}

\begin{abstract}
Isaak posed the following problem. Suppose $T$ is a tournament having a minimum feedback arc set which induces an acyclic digraph with a hamiltonian path. Is it true that the maximum number of arc-disjoint cycles in $T$ equals the cardinality of minimum feedback arc set of $T$? We prove that the answer to the problem is in the negative. Further, we study the number of arc-disjoint cycles through a vertex $v$ of the minimum out-degree in an oriented graph $D$. We prove that if $v$ is adjacent to all other vertices, then $v$ belongs to $\delta^+(D)$ arc-disjoint cycles.


\end{abstract}

\maketitle

\section{Introduction}
Let $D = (V, A)$ be a digraph.
A set of arcs $S \subseteq A$ is called a \textsl{feedback arc set} if $D-S$ is acyclic. The minimum number of elements in a feedback arc set of $D$ is denoted by $\tau(D)$. The maximum number of arc-disjoint cycles in $D$ is denoted by $\nu(D)$. Let $\pi = v_1,v_2,\ldots,v_n$ be an ordering of the vertices of $D$. An arc $v_iv_j \in A$ is called \textsl{backward with respect to $\pi$} if $i > j$; $\pi$ is \textsl{$\tau$-optimal} if the number of backward arcs with respect to $\pi$ is minimum among all orderings of vertices of $D$. For a vertex $v \in V$, we use the following notation:
$$
N^+(v) = \{u\in V-v : vu \in A \}, \quad  N^-(v) = \{u\in V-v : uv \in A \},
$$
$$
 N^{+2}(v) = \bigcup_{u \in {N^+(v)}} N^{+}(u) - N^+(v)-\{v\}.
$$
The \textsl{out-degree} of a vertex $v\in V$ is $d^+(v) = |N^+(v)|$, and  $\delta^+(D) = \min \{d^+(v): v \in V \}$. We use Bang-Jensen and Gutin \cite{bang} as reference for undefined terms. %For example: a circuit in a digraph is called a cycle.

It is well known that $S$ is a minimum feedback arc set in a digraph $D$ if and only if there exists a $\tau$-optimal ordering $\pi$ of vertices in $D$ such that $S$ is the set of backward arcs with respect to $\pi$ (see Bang-Jensen  and Gutin \cite{bang}). Hence it follows that for every digraph $D$ we have $\tau(D) \geq \frac{1}{2}\delta^+(D)(\delta ^{+}(D)+1)$ (see Remark 3).
Erd\H{o}s and Moon \cite{erdos} proved that for every $n \geq 3$ there exists a tournament $T_n$ of order $n$ such that $\tau(T_n) \geq \frac{1}{4}n(n-1)- \frac{1}{2}\sqrt {n^3\log_{e}n}$. A slightly better result was obtained by de la Vega in \cite{vega}.
% proved that for every $n \geq 10$ there exists a tournament $T_n$ of order $n$ such that $\tau(T_n) > [\frac{1}{3}n[\frac{1}{2}(n-1)]]$, answering a conjecture of Kotzig \cite{kotzig}.
On the other hand, it follows from a result by Chatrand, Geller and Hedetniemi \cite{chatrand} that $\nu(T_n) \leq \lfloor \frac{1}{3}n \lfloor\frac{1}{2}(n-1)\rfloor \rfloor$. Even though not always $\tau(D) = \nu(D)$, Isaak \cite{isaak}
conjectured (Conjecture 15.4.1 of \cite{bang}) that if $T$ is a tournament having a minimum feedback arc set which induces a transitive subtournament of $T$, then $\tau(T) = \nu(T)$. He posed also the following question (Problem 15.4.2 of \cite{bang}): Suppose $T$ is a tournament having a minimum feedback arc set which induces an acyclic digraph with a hamiltonian path. Is it true that $\tau(T) = \nu(T)$? We prove that the answer to the Isaak question is in the negative.

Mathematical sociologist Landau \cite{landau} proved that in every tournament $T$, if a vertex $v$ has the minimum out-degree, then it belongs to $\delta^+(T)$ different $3$-cycles. If $T$ is eulerian, then every vertex has the minimum out-degree. One may guess that for this case every vertex belongs to $\delta^+(T)$ arc-disjoint $3$-cycles. However, it is not true (see Remark 4). We prove (Theorem 2.1) that in every oriented graph $D$, if $v$ is a vertex of the minimum out-degree which is adjacent to all other vertices, then it belongs to $\delta^+(D)$ arc-disjoint cycles. Mader \cite{mader} consider linkages in digraphs with lower bounds on the out-degree. He proved a theorem which is related to our result: every digraph $D$ with $\delta^+(D) \geq n$ contains a pair of distinct vertices $x$, $y$ with $n-1$ arc-disjoint $(x,y)$-paths.

We conjecture that every tournament $T$ has a vertex of the minimum out-degree
which belongs to $\delta^+(T)$ arc-disjoint $3$-cycles. It is connected with the following conjecture put forward by Hoang and Reed \cite{hoang}: every digraph $D$ with $\delta^+(D) = n$ contains a sequence $C_1$, $C_2$, \ldots, $C_n$ of cycles such that $\bigcup^{j-1}_{i=1}C_i$ and $C_j$ have at most one vertex in common. In the case of $n = 2$ the Hoang and Reed conjecture was proved by Thomassen \cite{thomassen}. The conjecture was verified for tournaments by Havet, Thomass\'{e} and Yeo \cite{yeo}. They proved that every tournament $T$ with $\delta^+(T) = n$ contains a sequence $C_1$, $C_2$, \ldots, $C_n$ of $3$-cycles such that $\bigcup^{j-1}_{i=1}C_i$ and $C_j$ have exactly one vertex in common.  Our conjecture is also connected with the following conjecture posed by Seymour \cite{dean} (see also Seymour's Second Neighbourhood  Conjecture \cite{bang}): every oriented graph has a vertex $v$ such that $|N^+(v)| \leq |N^{+2}(v)|$. In the case of tournaments this conjecture was proved by Fisher \cite{fisher}. An elementary proof for the case of tournaments was found by Havet and Thomass\'{e} \cite{havet} (see also Bang-Jensen and Gutin  \cite{bang}).

\section{The Isaak problem}
Let $T$ be a tournament of order $13$ in Figure 1 and let $\alpha = a,b,c,d,e,f,g,h,i,j,k,l,m$ be an ordering of vertices of $T$. We will prove that $\nu(T) = 11$. Notice that the set of all backward arcs with respect to $\alpha$ is a feedback arc set which induces an acyclic digraph with a hamiltonian path $m, k, i, g, e, c, a$. We will prove that the ordering $\alpha$ is $\tau$-optimal. Hence, $\tau(T) = 12$

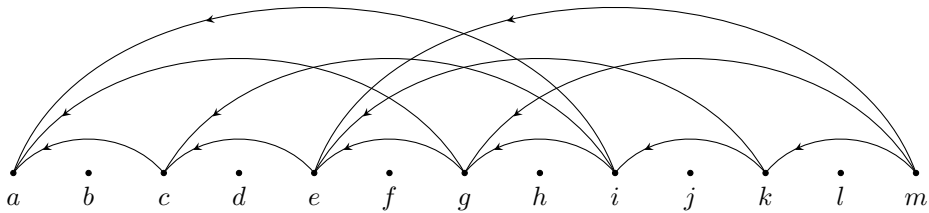
\begin{figure}[hb]%1
\centering
\begin{tikzpicture}[>=stealth, shorten >=1pt, bend angle=45]

\coordinate (a) at (0,0);
\coordinate (b) at (1,0);
\coordinate (c) at (2,0);
\coordinate (d) at (3,0);
\coordinate (e) at (4,0);
\coordinate (f) at (5,0);
\coordinate (g) at (6,0);
\coordinate (h) at (7,0);
\coordinate (i) at (8,0);
\coordinate (j) at (9,0);
\coordinate (k) at (10,0);
\coordinate (l) at (11,0);
\coordinate (m) at (12,0);

\foreach \i in {a,b,c,d,e,f,g,h,i,j,k,l,m}
\filldraw[fill=black] (\i) circle (1pt);
\foreach \i in {a,b,c,d,e,f,g,h,i,j,k,l,m}
\draw (\i) node[anchor=north]{$\strut\i$};

\begin{scope}[decoration={markings,
mark=at position -0.5cm with {\arrow[black]{stealth};}}]
\draw[postaction={decorate}] (c) to[in=50,out=130] (a);
\draw[postaction={decorate}] (e) to[in=50,out=130] (c);
\draw[postaction={decorate}] (g) to[in=50,out=130] (e);
\draw[postaction={decorate}] (i) to[in=50,out=130] (g);
\draw[postaction={decorate}] (k) to[in=50,out=130] (i);
\draw[postaction={decorate}] (m) to[in=50,out=130] (k);
\end{scope}

\begin{scope}[decoration={markings,
mark=at position -1cm with {\arrow[black]{stealth};}}]

\draw[postaction={decorate}] (g) to[in=60,out=120] (a);
\draw[postaction={decorate}] (i) to[in=60,out=120] (c);
\draw[postaction={decorate}] (k) to[in=60,out=120] (e);
\draw[postaction={decorate}] (m) to[in=60,out=120] (g);
\end{scope}

\begin{scope}[decoration={markings,
mark=at position .65 with {\arrow[black]{stealth};}}]

\draw[postaction={decorate}] (i) to[in=70,out=110] (a);
\draw[postaction={decorate}] (m) to[in=70,out=110] (e);
\end{scope}

\end{tikzpicture}

\caption{The tournament $T$. Only backward arcs with respect to the ordering $\alpha$ are shown.} %$\nu(T)= 11$ and $\tau(T) = 12$.}
\end{figure}

Let  $\mathcal{C}$ be a family of the following  arc-disjoint $3$-cycles in $T$:
 $$
 (a, b, c), (c, d, e), (e, f, g), (g, h, i), (i, j, k), (k, l, m), (a, d, g), (c, f, i), (e, h, k), (g,~j,~m), (a, e, i).
 $$
Suppose that $U$ is the union set of arcs of all cycles belonging to the family $\mathcal{C}$. Notice that $me$ is the only backward arc with respect to $\alpha$ which does not belong to $U$. If $me$ is an arc of some cycle in the tournament, then this cycle has an arc belonging to $U$. Hence, $\mathcal{C}$ is a maximal family of arc-disjoint cycles in $T$. It is easy to see that every backward arc determines uniquely a family of eleven arc-disjoint cycles omitting this backward arc. By analogy, we can check that every such family is a maximal family of arc-disjoint cycles in $T$. Hence,  $\mathcal{C}$ is a maximum family of arc-disjoint cycles in $T$. Thus $\nu(T) = 11$. Let  $S$ be a  minimum feedback arc set in $T$. If $S$ has eleven elements, then it satisfies the following conditions:
\begin{itemize}
\item [(1)]
every cycle in $\mathcal{C}$ has exactly one arc belonging to $S$,
\item [(2)]
$S \subset U$.
\end{itemize}
Since $im \notin U$ and $me \notin U$, by (2), the arc $ei$ of a cycle $(e, i, m)$ belongs to $S$. Hence, by (1), the arc $ia$ of a cycle $(a, e, i)\in \mathcal{C}$ does not belong to $S$. Since $ia \notin S$ and $ah \notin U$, by (2), the arc $hi$ of a cycle $(a, h, i)$ belongs to $S$. Hence, by (1), the arc $gh$ of a cycle $(g, h, i)\in \mathcal{C}$ does not belong to $S$. Since $gh \notin S$ and $hm \notin U$, by (2), the arc $mg$ of a cycle $(g, h, m)$ belongs to $S$. Hence, by (1), the arc $jm$ of a cycle $(g, j, m) \in \mathcal{C}$ does not belong to $S$. Thus we obtain a contradiction, a cycle $(e, j, m)$ is arc-disjoint with the set $S$. Hence, $S$ has at least twelve elements.

\begin{remark}
If we add three arcs $mc$, $kc$ and $ka$ to the tournament $T$ we obtain a tournament $T^{'}$ having a minimum feedback arc set which induces an acyclic digraph with a hamiltonian path, such that $\nu(T^{'})= 14$ and $\tau(T^{'}) = 15$.
\end{remark}

\begin{remark}
We checked case by case that if $T_{\leq6}$ is a tournament of order at most $6$, then the maximum number of arc-disjoint cycles in $T_{\leq6}$ is equal to the cardinality of a minimum feedback arc set of $T_{\leq6}$. We show in Figure 2 a tournament $T_7$ of order $7$ such that $\nu(T_7)= 4$ and $\tau(T_7) = 5$.
\end{remark}

\begin{figure}[h]%2
\centering
\begin{tikzpicture}[>=stealth, shorten >=1pt, bend angle=45]

\coordinate (a) at (0,0);
\coordinate (b) at (1,0);
\coordinate (c) at (2,0);
\coordinate (d) at (3,0);
\coordinate (e) at (4,0);
\coordinate (f) at (5,0);
\coordinate (g) at (6,0);

\foreach \i in {a,b,c,d,e,f,g}
\filldraw[fill=black] (\i) circle (1pt);
\foreach \i in {a,b,c,d,e,f,g}
\draw (\i) node[anchor=north]{$\strut\i$};

%\draw[->] (c) to[in=60,out=120] (a) [arrtwo];
\draw[->] (c) to[in=60,out=120] (a);
\draw[->] (e) to[in=60,out=120] (c);
\draw[->] (g) to[in=60,out=120] (d);
\draw[->] (f) to[in=60,out=120] (b);

%\pgfarrowsdeclarecombine[2cm]{artwo}{artwo}{latex}{latex}{}{stealth}

%\draw[-artwo ] (f) to[in=60,out=120]  (a);
%\draw (f) to[in=60,out=120]  (a) node[pos=0.5] (fa) {} ;
%\draw[->] (fa) to[in=0,out=180]  (fa);

\begin{scope}[decoration={markings,
%mark=at position 1cm with {\node[red]{1cm};},
%mark=at position .75 with {\arrow[blue,line width=2mm]{>}},
mark=at position -1cm with {\arrow[black]{stealth};}}]
\draw[postaction={decorate}] (f) to[in=60,out=120] (a);
\end{scope}

\end{tikzpicture}
\caption{The tournament $T_7$. Only backward arcs with respect to an ordering $a,b,c,d,e,f,g$ are presented.} %$\nu(T_7)= 4$ and $\tau(T_7) = 5$.}

\end{figure}
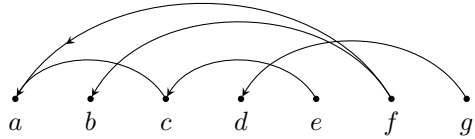

\begin{remark}
If $D$ is a digraph, then $\tau(D) \geq \frac{1}{2}\delta^+(D)(\delta ^{+}(D)+1)$. We proceed by induction on the order of $D$. Let $S$ be a feedback arc set in $D$ of minimum size, and suppose that $\pi = v_1, v_2,\ldots, v_n$ is a $\tau$-optimal ordering of vertices in $D$ such that $S$ is the set of backward arcs with respect to $\pi$. Since $v_1, v_2,\ldots, v_{n-1}$  is a $\tau$-optimal ordering of vertices in $D-v_n$, $S' = S -\{v_nv: v \in N^+(v_n)\}$ is a feedback arc set in $D-v_n$ of minimum size. Notice that $\delta^+(D-v_n) \geq \delta^+(D)-1$. Hence, the size of $S'$ is at least $\frac{1}{2}(\delta^+(D)-1)\delta ^{+}(D)$. Therefore, the size of $S$ is at least
$$\frac{1}{2}(\delta^+(D)-1)\delta ^{+}(D) + \delta ^{+}(D)= \frac{1}{2}\delta^+(D)(\delta ^{+}(D)+1).$$
\end{remark}

\section{Arc-disjoint cycles through a vertex of the minimum out-degree}
For a pair $X$, $Y$ of vertex sets of an oriented graph $D = (V, A)$, we define
$$(X,Y)_D = \{(x,y)\in A : x \in X, y \in Y\}.$$
\begin{theorem}\label{theorem 2.1}
Let $D$ be an oriented graph, and suppose that $v_0$ is a vertex which is adjacent to all other vertices in $D$. Let $a = \min \{d^+(v) : v\in N^+(v_0)\}$ and  $b = \min \{d^+(v) : v\in N^-(v_0)\}$.
If $d^+(v_0) \leq \min (a, \frac{1}{2}(a+b+1))$, then $v_0$ belongs to $d^+(v_0)$ arc-disjoint cycles.
\end{theorem}
\begin{PrfFact}
Let $\Gamma = \{\gamma_1, \gamma_2, \ldots, \gamma_m\}$ be a maximum family of arc-disjoin cycles through $v_0$.  Let
$\gamma_i = v_0 v^{i}_{1}\ldots v^{i}_{n(i)} v_0$, for $i = 1, \ldots, m$. By Menger's theorem (see \cite{menger} and \cite{bang} there exists a set $\Delta$ of $m$ arcs covering all cycles containing the vertex $v_0$. Suppose that $k$ is the number of arcs in  $\Delta$ with the head $v_0$. If $k > 0$, we can assume that $(v^{i}_{n(i)},v_0) \in \Delta$ if and only if $1 \leq i \leq k$. Let us denote $K = \{v^1_1, v^2_1,\ldots, v^k_1 \}$, $X = \{v^1_{n(1)}, v^2_{n(2)},\ldots, v^k_{n(k)} \}$ (if $k = 0$ we set $K = X = \emptyset$), $L = \{v^{k+1}_1, v^{k+2}_1,\ldots, v^{m}_1 \}$, $Y =  N^-(v_0)-X$, and $M = N^+(v_0)-K-L$. First we prove the following inequality:
\begin{itemize}
\item [(1)]
$|(K\cup X \cup M, Y)_D| \leq |(L, K\cup X \cup M)_D|$.
\end{itemize}
Assume that $(v,y) \in (K\cup X \cup M, Y)_D$. Notice that $(y, v_0) \notin \Delta$. If $v \in K \cup M$, then $(v_0, v) \notin \Delta$. Hence, the arc $(v,y)$ of a cycle $v_0vyv_0$ belongs to  $\Delta$. If $v \in X$, then $v = v^i_{n(i)}$, for some $i \leq k$. Hence, the arc $(v^i_{n(i)}, y)$ of a cycle $v_0 v^{i}_{1}\ldots v^{i}_{n(i)}yv_0$ belongs to $\Delta$. Thus, $(v,y)$ is an arc of some cycle $\gamma_i$, for $i > k$. Accordingly, to every arc $(v,y) \in (K\cup X \cup M, Y)_D \subseteq \Delta$ we can assign an arc $(l,v') \in (L, K\cup X \cup M)_D$ such that $(l, v')$ and $(v,y)$ belong to the same cycle in $\Gamma$. The above assignment is injective, because different arcs in $\Delta$ belong to arc-disjoint cycles in $\Gamma$. Hence, (1) holds.

Let us complete the oriented graph $D$ to a tournament $T$ with the same vertex set $V$. Hence, $|K\cup X \cup M, V-Y)_T| \geq  |(K\cup X \cup M, V-Y)_D|$. Therefore, by (1) we obtain:
\begin{align*}
|(L, K\cup X \cup M)_D | & + |(K\cup X \cup M, V)_T|
\geq |K\cup X \cup M, Y)_D| + |(K\cup X \cup M, V)_T|
\cr &
\geq |(K \cup X \cup M, Y)_T| + |(K \cup X \cup M, V)_D|.
\end{align*}
Accordingly,
\begin{align*}
|K \cup X \cup M|(|V|-1) & = |(V, K \cup X \cup M)_T| + |(K \cup X \cup M, V)_T|
\cr &
= |(V-L, K \cup X \cup M)_T| + |(L, K \cup X \cup M)_T| +|(K \cup X \cup M, V)_T| \cr
& \geq  |( V-L, K \cup X \cup M)_T| +|(K \cup X \cup M, Y)_T| +|(K \cup X \cup M, V)_D| \cr
&
= |(K \cup X \cup M, K \cup X \cup M )_T| +  |(\{v_0\}, K \cup X \cup M)_T|  \cr
& \quad + |(Y, K \cup X \cup M)_T| + |(K \cup X \cup M, Y)_T|
  +  |(K \cup X \cup M, V)_D|
\cr & \geq \textstyle\frac{1}{2}|K \cup X \cup M|(|K \cup X \cup M|-1) + |K| + |M|  \cr
& \quad
+ |K \cup X \cup M||Y|+  a|K| + b |X| + a|M|.
\end{align*}
Since $|V| - 1 = d^+(v_0) + |X| + |Y|$ and $|K| = |X|$, we have
\begin{align*}
2|K|d^+(v_0) + |M|d^+(v_0) \geq \left(|K|+ \textstyle\frac{1}{2}|M|\right)(|M|-1) + |M| + (a+b+1)|K| + a|M|.
\end{align*}
Thus, $|M| = 0$, because $(a+b+1) \geq 2d^+(v_0)$ and $a \geq d^+(v_0)$.
\end{PrfFact}

\begin{remark}

\begin{figure}%3
\centering
\begin{tikzpicture}[>=stealth, shorten >=1pt, bend angle=45]
\coordinate (a) at (0,0);
\coordinate (b) at (1,0);
\coordinate (c) at (2,0);
\coordinate (d) at (3,0);
\coordinate (e) at (4,0);
\coordinate (f) at (5,0);
\coordinate (g) at (6,0);
\coordinate (h) at (7,0);
\coordinate (i) at (8,0);
\coordinate (j) at (9,0);
\coordinate (k) at (10,0);

\foreach \i in {a,b,c,d,e,f,g,h,i,j,k}
\filldraw[fill=black] (\i) circle (1pt);
\foreach \i in {a,b,c,d,e,f,g,h,i,j,k}
\draw (\i) node[anchor=north]{$\strut\i$};

\draw[->] (c) to[in=45,out=135] (a);
\draw[->] (g) to[in=45,out=135] (d);
\draw[->] (j) to[in=45,out=135] (h);

\draw (-0.2, -0.1) rectangle (2.4,0.7);
\draw (6.8, -0.1) rectangle ++(2.4,0.8);
\node (ac) at (1,0.6) {};
\node (hj) at (8,0.58) {};
\draw[->, shorten >=-1pt] (hj) to[in=45,out=135] (ac);

\draw (-0.4, -0.15) rectangle (4.4,1.0);
\node (ae) at (0,0.85) {};
\draw[->, shorten >=1pt] (k) to[in=45,out=110] (ae);

\end{tikzpicture}

\caption{The  eulerian tournament $T_{11}$.\newline $(\{h, i, j\}, \{a,b,c\})_{T_{11}} \cup (\{k\}, \{a,b,c,d,e\})_{T_{11}} \cup \{(c,a), (g,d), (j,h)\}$ is the set of all backward arcs with respect to the ordering $\beta$.}
\end{figure}
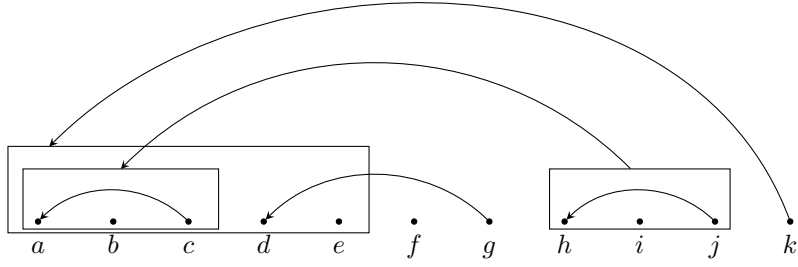

 Let $T_{11}$ be the eulerian tournament in Fig.~3, and $\beta = a,b,c,d,e,f,g,h,i,j,k$ an ordering of its vertices. Suppose that
$$(\{h, i, j\}, \{a,b,c\})_{T_{11}} \cup (\{k\}, \{a,b,c,d,e\})_{T_{11}} \cup \{(c,a), (g,d), (j,h)\}$$
is the set of all backward arcs with respect to the ordering $\beta$. Let $k v_1 v_2$ be a $3$-cycle through the vertex $k$. Notice that, if $v_1 \in \{a, b,c\}$, then $v_2 \in \{f, g\}$. Hence, the vertex $k$ does not belong to $\delta^+(T_{11})$ arc-disjoint $3$-cycles.
\end{remark}\


\begin{thebibliography}{99}

\bibitem{bang} J. Bang-Jensen and G. Z. Gutin, Digraphs, Springer Monographs in Mathematics, Second Edition (2010).

\bibitem{chatrand} G. Chatrand, D. Geller and S. Hedetniemi, Graphs with forbidden subgraphs, J. Combin. Theory Ser. B 10 (1971), 12--41.

\bibitem{dean} N. Dean and B. J. Latka, Squaring the tournaments --- an open problem,
Congr. Numer. 109 (1995), 73--80.

\bibitem{erdos} P. Erd\H{o}s and J. W. Moon, On sets of conistent arcs in tournament, Can. Math. Bull. 8 (1965), 269--271.

\bibitem{fisher} D. C. Fisher, Squaring a tournament: a proof of Dean's conjecture,
J. Graph Theory 23(1) (1996), 43--48.

\bibitem{havet} F. Havet and S. Thomass\'{e}, Median orders of tournaments: a tool for the second neighbourhood problem and Sumner's conjecture, J. Graph Theory 35(4)     (2000), 244--256.

\bibitem{yeo} F. Havet, S. Thomass\'{e} and A. Yeo, Ho\`{a}ng-Reed Conjecture holds for tournaments. Discrete Math. 308 (2008), 3412--3415.

\bibitem{hoang} C. T. Ho\`{a}ng and B. Read, A note on short cycles in digraphs, Discrete Math. 66 (1--2) (1987), 103--107.

\bibitem{isaak} G. Isaak, Tournaments as feedback arc sets, Electron J. Combin. 2 (1995), 1--19.

\bibitem{landau} H. G. Landau, On dominance relations and the structure of an animal
 societies III. The condition for a~score structure, Bull. Math. Biophys. 15 (1953), 143--148.

\bibitem{mader} W. Mader, Degree and local connectivity in digraphs, Combinatorica
     5(2) (1985), 161-165.

\bibitem{menger} K. Menger, Zur allgemeinen Kurventheorie, Fund Math. 10 (1927), 96--115.
%\bibitem{sullivan} B. Sullivan, A summary of results and problems related to the
%Cacceta-H\"{a}ggkvist conjecture. Preprint arXiv:math.C0/0605646v1, (May 2006).

\bibitem{thomassen} C. Thomassen, The 2-linkage problem for acyclic digraphs, Discrete Math. 55 (1) (1985), 73--87.

%\bibitem{bermond} J. C. Bermond, The circuit-hypergraph of a tounament, Infinite and Finite sets, Vol.I (ed. A. Hajnal et al. ), North-Holand, Amsredam, (1975), 165--180.


\bibitem{vega}W. F. de la Vega, On the maximum cardinality of a consistent set of arcs in a~random tournament, J.~Combin. Theory Ser. B 35 (1983), 328--332.

%\bibitem{barthelemy} J. P. Barth\'{e}l\'{e}my, O. Hudry, G. Isaak, F.S. Roberts, and B. Tesman. The reversing number of a digraph. Discrete Appl. Math. 60(1-3), (1995), 39--76.

%\bibitem{kotzig} A. Kotzig, On the number of cyclicity of antisymmetric directed graphs, Discrete Math. 12, (1975), 17--25.


\end{thebibliography}
\end{document}